\documentclass[11pt]{amsart}
\usepackage{amssymb,amscd,amsthm,amsmath,amsfonts}

\newcommand{\norm}[1]{\mbox{$\left\|#1\right\|$}}
\newcommand{\x}{\times}
\newcommand{\cs}{\mbox{$C^{*}$-algebra}}
\newcommand{\css}{\mbox{$C^{*}$-algebras}}

\newcommand{\Z}{\mathbb{Z}}
\newcommand{\Pos}{\mathbb{P}}
\newcommand{\ov}[1]{\mbox{$\overline{#1}$}}

\newcommand{\al}{\mbox{$\alpha$}}

\newcommand{\bt}{\mbox{$\beta$}}
\newcommand{\ga}{\mbox{$\gamma$}}

\newcommand{\de}{\mbox{$\delta$}}

\newcommand{\la}{\mbox{$\lambda$}}

\newcommand{\La}{\mbox{$\Lambda$}}

\newcommand{\si}{\mbox{$\sigma$}}

\newcommand{\mfD}{\mathfrak{D}}

\newcommand{\mfX}{\mathfrak{X}}

\newcommand{\mfR}{\mathfrak{R}}
\newcommand{\mfH}{\mathfrak{H}}

\newcommand{\calE}{\mathcal{E}}

\newcommand{\bgc}{\begin{center}}

\newcommand{\edc}{\end{center}}
\newcommand{\be}{\begin{enumerate}}
\newcommand{\ee}{\end{enumerate}}
\newcommand{\beqn}{\begin{eqnarray}}
\newcommand{\eeqn}{\end{eqnarray}}
\newcommand{\beqns}{\begin{eqnarray*}}
\newcommand{\eeqns}{\end{eqnarray*}}
\newcommand{\bq}{\begin{quote}}
\newcommand{\eq}{\end{quote}}

\newcommand{\bi}{\begin{itemize}}
\newcommand{\ei}{\end{itemize}}
\newcommand{\bd}{\begin{description}}
\newcommand{\ed}{\end{description}}
\newcommand{\ra}{\mbox{$\rightarrow$}}

\theoremstyle{plain}
\newtheorem{theorem}{Theorem}

\newtheorem{definition}{Definition}
\newtheorem{proposition}{Proposition}
\newtheorem{corollary}{Corollary}

\numberwithin{equation}{section}

\begin{document}
\title{Graph inverse semigroups, groupoids and their $C^{*}$-algebras}
\author{Alan L. T. Paterson}
\address{Department of Mathematics\\
University of Mississippi\\
University, MS 38677}
\email{mmap@olemiss.edu}
\date{}
\begin{abstract}
There is now a substantial literature on graph $\css$. Under a locally finite
condition on a countable,
directed graph, Kumjian, Pask, Raeburn, Renault
showed that the $\cs$ of the graph can be realized as the $\cs$ of the {\em
path groupoid}, i.e. the groupoid determined by the infinite paths in the
graph. In the present paper, we remove the local finiteness requirement. The
path groupoid in the general context is obtained through the universal groupoid
of a certain inverse semigroup associated with the graph. This inverse
semigroup is called the {\em graph inverse semigroup}, and graph
representations turn out to be just representations of this inverse semigroup.
A certain reduction of the universal groupoid gives the {\em path groupoid} of
the graph, and its $\cs$ is isomorphic to the $\cs$ of the graph. The unit
space of the path groupoid contains the infinite paths of the graph, but also
contains some finite paths. We show that, as in the locally finite case, the
path groupoid is always amenable, and we give a groupoid proof of a recent
theorem of W. Szymanski, characterizing when a graph $\cs$ is simple.
\end{abstract}

\maketitle

\section{Introduction}
There is now a substantial literature on the \cs\ associated with a (countable)
directed
graph. Among the motivations for this are the Bratteli diagrams for
AF-algebras and especially the Cuntz-Krieger algebras, which can be interpreted
in graphical terms with the matrix $A$, determining such an algebra, regarded as
the incidence matrix of the graph. A helpful survey of the field up to 1998 has
been given by A. Kumjian (\cite{Kumjian}). In the literature surveyed by
Kumjian, the graphs considered satisfy finiteness conditions on the edges
coming into and going out of a vertex. Usually, {\em local finiteness} has been
required, viz. that at each vertex $v$, only finitely many edges end at $v$,
and only finitely many edges start at $v$. Many of these results were
shown to be valid in the {\em row finite} case, i.e. where we only assume that
finitely many edges start at $v$ for each vertex $v$ of the graph.

Recently, there has been interest in removing the finiteness requirements on
the graph. See, for example,
\cite{BPRS,EL,EL2,FLR,HR,RS,ren,Szy1,Szy2,Szy3,Szy4}.
Motivation for this
is to be found even in the original paper \cite{CuntzSimple} of Cuntz in which
he studied what are now called the Cuntz algebras $\mathcal{O}_{n}$. These are
graph $\css$, the natural graph associated with $\mathcal{O}_{n}$ being a
bouquet of $n$ circles (one vertex). When $n<\infty$, the graph is locally
finite. When $n=\infty$, this is no longer the case. (Indeed up to the work
of Exel and Laca (\cite{EL}), $\mathcal{O}_{\infty}$ was the only
known example that had been worked out of a Cuntz-Krieger algebra whose
matrix was not row-finite.) Even in the case of the Cuntz algebras, there is a
significant difference between the row finite/non-row finite cases.
This is shown by the fact that in \cite{CuntzSimple}, in a number of instances,
proofs given for $\mathcal{O}_{n}$ with $n<\infty$ do not apply in the case
of $\mathcal{O}_{\infty}$, so that different proofs were required for that case.

A representation of a directed graph $\calE$ is determined by certain partial
isometries on a Hilbert space indexed by the vertices and edges. The $\cs$
$C^{*}(\calE)$ of $\calE$ is the universal $\cs$ for such representations.
In the locally finite case, insight into this $\cs$ is provided by the {\em
path groupoid} $G$ of $\calE$. This is an r-discrete groupoid, and its units
are just the infinite paths in $\calE$. Then (cf. \cite{KPRR}) if $\calE$ is
locally (or even row) finite and has no sinks, we have $C^{*}(\calE)\cong
C^{*}(G)$, and this groupoid connection ennables one to interpret properties
of $C^{*}(\calE)$ in dynamical terms.

When we look for this groupoid connection without assuming local finiteness, we
run into a problem since the space of infinite paths in $\calE$ is not a
locally compact Hausdorff space in any obvious way. The main objective of
this paper is to show that, despite this problem with the infinite paths, there
is a natural path groupoid available in complete generality.  A peculiarity
of this groupoid is that we have to include some {\em finite} paths.  To
determine what this groupoid is, we use an approach based on inverse semigroups.

Indeed, there is a close connection between inverse semigroups and r-discrete
groupoids. This was first observed and studied by J. Renault in his monograph
\cite{rg}. Recently, the present author in his book \cite{Paterson} gave a
detailed account of this connection. One of the main results of
\cite[Ch.4]{Paterson} is a construction which associates with any (countable)
inverse semigroup $S$ a {\em universal} groupoid $H$ canonically related to $S$
such that $C^{*}(S)\cong C^{*}(H)$. The path groupoid will be a reduction of
the universal groupoid for a certain inverse semigroup associated with the
graph. So we have to come up with an appropriate inverse semigroup $S$ for
$\calE$.

In the $\mathcal{O}_{n}$-case, the appropriate inverse semigroup is $S_{n}$
which was introduced by Renault in \cite[p.141]{rg}. For a general directed
graph $\calE$, a natural candidate for such an inverse semigroup is the inverse
semigroup of path pairs $(\al,\bt)$ considered in \cite[p.158]{Paterson}. As
mentioned there, an unsatisfactory feature of this inverse semigroup is that,
unlike the Cuntz case, it is not generated by the edges. The way out of this
difficulty is suggested by the definition of a representation of $\calE$
referred to above. {\em We have to include the vertices in the inverse
semigroup.}

The precise details of this are given in \S2. The {\em graph inverse semigroup}
of $\calE$, denoted by $S_{\calE}$, is defined as the inverse semigroup
generated by the edges and vertices of $\calE$ subject to certain
``Cuntz-like'' relations. We show that $S=S_{\calE}$ can be identified with an
inverse semigroup of path pairs, where we permit vertices to stand as
surrogates for ``paths'' of zero length. This ``path-pair'' characterization
of $S$ is useful for calculation purposes.

The universal groupoid $H$ of $S$ is then determined in the paper. We
focus on a natural reduction of this which we will call the {\em graph
groupoid} $G$. This groupoid is similar to the graph groupoid in the locally
finite case and indeed is the same groupoid in that case. The main difference
is that we must include in the unit space $G^{0}$ those {\em finite} paths
ending in a vertex $v\in V_{\infty}$, where $V_{\infty}$ is the set of vertices
emitting infinitely many edges. Sequences of infinite paths can converge to
such finite paths, in particular to $V_{\infty}$-vertices.   We write $V_{f}$
for $V\setminus V_{\infty}$, the set of vertices emitting only finitely many
edges.

The representations of $\calE$ correspond to a certain class of representations
of $S$, and using the canonical identification of $C^{*}(S)$ with $C^{*}(H)$
in \cite[Ch. 4]{Paterson}, we show that in fact these representations
correspond to the representations of $G$. So we obtain that
$C^{*}(\calE)=C^{*}(G)$ just as in the locally finite case.

Next we turn to the issue of the amenability of $C^{*}(\calE)$.  A groupoid
form of Cuntz's crossed-product argument for the nuclearity (and hence
amenability) of $\mathcal{O}_{n}$ was given by J. Renault in his book
\cite[p.138ff.]{rg}. For the case of a locally finite graph $\calE$, an
ingenious adaptation of this groupoid argument was given by Kumjian, Pask,
Raeburn and Renault (\cite{KPRR}) to show that $C^{*}(\calE)$ is always
amenable.  We adapt this argument further to establish the same result in
complete generality. The technical difficulty is that of dealing with the paths
of finite length in $G^{0}$.

The last section of the paper deals with the {\em simplicity} of
$C^{*}(\calE)$. Kumjian, Pask, Raeburn and Renault (\cite{KPRR}) characterized
the simplicity of $C^{*}(\calE)$ in terms of {\em condition (K)} and {\em
cofinality} in the locally finite case. (See below.)
Recently, W. Szymanski (\cite[Theorem 12]{Szy1})
characterized the simplicity of $C^{*}(\calE)$ for general $\calE$ using
the approach of Exel and Laca (\cite{EL}).
He showed that $C^{*}(\calE)$ is simple
if and only if ($\al$) and ($\bt$) hold where:
\bi
\item[($\al$)] all loops in $\calE$ have exits;
\item[($\bt$)] for any $v\in V$, the smallest hereditary, saturated subset of
$V$ that contains $v$ is equal to $V$.
\ei
(In $(\bt)$, a subset $H$ of $V$ is called {\em hereditary} if $w\in H$
whenever there is a path from some $v\in H$ to $w$, and is called {\em
saturated} if $w\in H$ whenever $w\in V_{f}$ and every edge starting at $w$
ends in $H$.)

We show first that for general $\calE$, the condition (K) of \cite{KPRR} -
that there are no vertices emitting exactly one loop - is equivalent
to $G$ being essentially
principal (i.e. every closed, invariant subset of $G^{0}$ contains a
dense subset of elements with trivial isotropy).
(This result is a counterpart to the result of Kumjian, Pask and
Raeburn (\cite[Lemma 3.4]{KPR}) that, in the locally finite case (no sinks),
the groupoid $G$ is essentially free if and only if every loop in $\calE$
has an exit.) We then show (Theorem~\ref{th:simple})
that $C^{*}(\calE)$ is simple if (a), (b) and (c)
hold, where:
\bi
\item[(a)] $\calE$ has property (K);
\item[(b)] $\calE$ is cofinal (i.e. given any $v\in V$ and any infinite path
$z$, then there is a path from $v$ to at least one of the vertices that $z$
passes through);
\item[(c)] if $v\in V$ and $w\in V_{\infty}$, then there is a path from $v$ to
$w$.
\ei
It is not difficult to show that (a),(b) and (c) together are equivalent to
($\al$) and ($\bt$) together, so that Theorem~\ref{th:simple} is effectively
just one of the implications in Szymanski's theorem.  The point of including
Theorem~\ref{th:simple} in the present paper is that it gives an
alternative approach to the simplicity of $C^{*}(\calE)$ using the groupoid
$G$ of the paper. The groupoid proof that we give is just a variant of the
method used by J. Renault in his book \cite{rg} (see also the
approach of \cite{KPRR}) to establish the simplicity of
the Cuntz $\css$:
using $G$ in place of the Cuntz groupoids, Renault's approach
works for a completely general graph $\calE$.

Another approach to the simplicity of graph $\css$ has been given by J\"{u}rgen Schweizer (\cite{Schweiz}). Schweizer introduces the concept of a {\em continuous measured diagram}, in which we are given locally compact spaces $D,X,Y$ with continuous maps $s:D\to X$, $r:D\to Y$, together with certain measures on $D$
indexed by the elements of $Y$. The case where $X=Y$ gives a natural $C^{*}$-correspondence, which in turn determines a {\em Cuntz-Pimsner} 
algebra $\mathcal{O}_{\mfD}$.  The discrete case of such a diagram is just that of a 
{\em directed graph} with set of vertices $X$, and in that case, $\mathcal{O}_{\mfD}$ 
is the same as the graph $\cs$.  Using \cite{Schweiz2}, 
Schweizer (\cite[Corollary 5.2]{Schweiz}) proves that when $X$ is compact and $r$ is onto, 
$\mathcal{O}_{\mfD}$ is simple if and only if the  
continuous measured diagram is non-periodic and minimal.  The simplicity results of Schweizer and Szymanski are different.  On the one hand, Schweizer's result is not restricted to the discrete (graphical) case, while on the other hand, its  compactness condition on $X$ in the graphical situation forces the vertex set to be finite. It would be interesting to relate the approach of Schweizer to the groupoid approach.

The author is grateful to Iain Raeburn and David Pask for bringing to his
attention the theorem of Wojciech Szymanski described above, and for pointing
out an error in an earlier version of the paper. He is also grateful to Marcelo
Laca, David Pask, Iain Raeburn, J\"{u}rgen Schweizer and Wojciech Szymanski for sending him recent
papers on graph $\css$. Finally, he wishes to thank Alan Donsig, John Meakin, David Pitts 
and Baruch Solel for helpful conversations on the topics of this paper.


\section{Some preliminaries}
Let $\calE$ be a directed graph. Its set of vertices is denoted by $V$ and its
set of edges by $E$. For each edge $e$, $s(e)$ is the initial vertex (source)
of $e$ and $r(e)$ the terminal vertex of $e$. It is assumed throughout the
paper that $V$ and $E$ are countable sets.

The graph $\calE$ is called {\em row finite} if for each $v\in V$, the set of
edges starting at $v$ is finite.  The graph is called {\em locally finite}
if for each $v\in V$, the set of edges starting at $v$ is finite, and the set
of edges terminating at $v$ is finite.  A vertex $v$ is called a {\em sink}
if there are no edges starting at $v$.  Throughout the paper we will assume
that there are no sinks in $\calE$ (unless the contrary is explicitly
specified.)

A reason for this is as follows. We want to examine $C^{*}(\calE)$ using a
``path groupoid'' $G$, and the paths that make up the unit space of the
appropriate groupoid $G$ are such that they cannot end at a sink. So if we want
to use the path groupoid to study $C^{*}(\calE)$, then we cannot allow sinks
in $\calE$. The $\cs$ $C^{*}(\calE)$ still makes good sense whether there are
sinks or not, and indeed, groupoid techniques can still be used using the {\em
universal groupoid} $H$ introduced in the paper rather than the path groupoid
$G$. Further, as is pointed out in \cite{BPRS}, it is technically easy to
reduce to the case where there are no sinks.

The set $V$ of vertices is the disjoint union of the sets $V_{f}$ and
$V_{\infty}$, where $v\in V_{f}$ if there are finitely many edges starting
at $v$, while $v\in V_{\infty}$ if there are infinitely many edges starting
at $v$.

The edge $e$ with reversed orientation is denoted by $e^{*}$. So
$s(e^{*})=r(e)$ and $r(e^{*})=s(e)$. We write $E^{*}$ as the set of $e^{*}$'s.

A {\em finite path} is a sequence $\al$ of edges $e_{1},\ldots ,e_{k}$ where
$s(e_{i+1})=r(e_{i})$ for $1\leq i\leq k-1$. We write
$\al=e_{1}e_{2}\ldots e_{k}$. The length $l(\al)$ of $\al$ is
just $k$. Each vertex $v$ is regarded as a finite path of length $0$. We
define $r(\al)=r(e_{k})$ and $s(\al)=s(e_{1})$. For $v\in V$, we set
$r(v)=v=s(v)$. The set of finite paths in $\calE$ is denoted by $Y$. Note that
$V\subset Y$. The set of infinite paths $z=z_{1}z_{2}\ldots $ in $\calE$
is denoted by $Z$. The length $l(z)$ of $z\in Z$ is defined to be
$\infty$.

Given two vertices $v,w$, we say that $v\geq w$ if there exists a
path $\al$ such that $s(\al)=v, r(\al)=w$. If $z$ is a path in $\calE$ and
$v\in V$, then we write $z\geq v$ to mean that $r(z_{i})\geq v$ for all $i$. We
write $z\leq v$ if $r(z_{i})\leq v$ for all $i$.

We now briefly discuss {\em inverse semigroups}.  A semigroup $S$ is called
an inverse semigroup if for all $s\in S$, there exists a unique element
$t\in S$ such that $sts=s,tst=t$.  We write the element $t$ as $s^{*}$.
Note that $s^{**}=s$. Inverse semigroups can be identified with semigroups of
partial one-to-one maps on a set that contain their inverses. Every element
$ss^{*}$ belongs to the set $I$ of idempotents of $S$. The set $I$ is a
commutative subsemigroup of $S$ and so is a semilattice. (In semigroup
theory, the semilattice of idempotents in $S$ is usually denoted by $E$, but
in this paper, $E$ stands for the set of edges in $\calE$.) Homomorphic images
of inverse semigroups are inverse semigroups.

A {\em congruence} $\rho$ on a semigroup $S$ is an equivalence relation
for which $(xa,xb),(ax,bx)\in \rho$
whenever $(a,b)\in \rho$.  If $\rho$ is a
congruence on $S$, then $S/\rho$ is a semigroup in the natural way.

We now recall the definition (\cite[pp.40-41]{CP1}) of the semigroup
$S$ generated by a set $X$ subject to the generating relations
$u_{\la}=v_{\la}$ ($\la\in \La$). Let $F_{X}$ be the semigroup of words in
$X$ with juxtaposition as product, and for $\la\in \La$, let
$u_{\la}, v_{\la}\in F_{X}$.
Then $S$ is just the quotient $F_{X}/\rho$ where $\rho$ is the smallest
congruence on $F_{X}$ containing all of the pairs $(u_{\la}, v_{\la})$.

Throughout the paper, we will assume the theory of r-discrete groupoids as
presented, for example, in \cite{rg,Paterson}.

\section{The inverse semigroup of a graph and the path groupoid}
Let $z$ be an element not in $V\cup E\cup E^{*}$. The {\em graph inverse
semigroup} $S=S_{\calE}$ of the graph $\calE$ is defined to be the semigroup
generated by $(V\cup E\cup E^{*})\cup\{z\}$ subject to the relations:
\bi
\item[(i)] $z$ is a zero for $S$;
\item[(ii)] $s(e)e=e=er(e), r(e)e^{*}=e^{*}=e^{*}s(e)$ for all $e\in E$;
\item[(iii)] $ab=z$ if $a,b\in V\cup E\cup E^{*}$ and $r(a)\neq s(b)$;
\item[(iv)] $e^{*}f=z$ if $e,f\in E$ and $f\neq e$.
\ei
Note that $S$ is countable since $E\cup V$ is.

We will see in a moment that $S$ is actually an inverse semigroup.  But for
the present, we describe a model $T$ for $S$.
The elements of $T$ are the pairs of finite paths $(\al,\bt)$ with
$r(\al)=r(\bt)$, together with a zero element $z$.  An involution is defined
on $T$ by: $z^{*}=z, (\al,\bt)^{*}=(\bt,\al)$, and the product is determined by:
\beqn
 (\al,\al'\mu)(\al',\bt) & = & (\al,\bt\mu)  \label{eq:prod1} \\
(\al,\al')(\al'\mu,\bt') & = & (\al\mu,\bt') \label{eq:prod2}
\eeqn
All other products in $T$ are $z$. Note that for finite paths $\al,\bt$, the
expression $\al\bt$ is defined only when $r(\al)=s(\bt)$, and in that case, is
just the finite path obtained in the obvious way by juxtaposition. Recall
that each vertex $v$ is regarded as a finite path. For any
finite path $\mu$, we specify that $s(\mu)\mu=\mu=\mu r(\mu)$.
Note that the product of two pairs $(\al_{1},\bt_{1}), (\al_{2},\bt_{2})$ is
non-zero if and only if $\bt_{1}$ is an initial segment of $\al_{2}$ or
conversely. In this case, we say that $\bt_{1}$, $\al_{2}$ are {\em
comparable}. Note that if $\al\mu, \bt\nu$ are comparable,
then so are $\al,\bt$. The set $I$ of
idempotents of $S$ is just the set of pairs $(\al,\al)$, and so can be
identified with $Y$. The product on $I$ of two elements $\al,\bt\neq z$ is
given by $\al\bt=\al=\bt\al$ if $\bt$ is an initial segment of $\al$, and is
$z$ if $\al,\bt$ are not comparable.

A similar semigroup to $T$ was suggested by the author in
\cite[p.158]{Paterson} but there is a significant difference: in
\cite{Paterson}, the paths $\al,\bt$ were of finite, non-zero length, and
the edges and vertices do not lie in the semigroup. However, as we will
see (Proposition~\ref{prop:Phi}),
the edges and vertices lie inside and generate $T$.

\begin{proposition}        \label{prop:T}
The set $T$, with the above involution and product, is an inverse semigroup.
\end{proposition}
\begin{proof}
We first have to show that the product in $T$ is associative. If one of the
terms is $z$, then this is obvious. So let
$(\al_{i},\bt_{i})\in T$ ($1\leq i\leq 3$),
$w_{1}=(\al_{1},\bt_{1})((\al_{2},\bt_{2})(\al_{3},\bt_{3}))$ and
$w_{2}= ((\al_{1},\bt_{1})(\al_{2},\bt_{2}))(\al_{3},\bt_{3})$.  We have to
show that $w_{1}=w_{2}$.

The cases that give $w_{1}\neq z$ are as follows (for appropriate words
$\mu,\nu$):
\bi
\item[(1)] $\bt_{2}=\al_{3}\mu$, $\bt_{1}=\al_{2}\nu$,
$w_{1}=(\al_{1},\bt_{3}\mu\nu)$;
\item[(2)]  $\bt_{2}=\al_{3}\mu$, $\al_{2}=\bt_{1}\nu$,
$w_{1}=(\al_{1}\nu,\bt_{3})$;
\item[(3)] $\al_{3}=\bt_{2}\mu$, $\bt_{1}=\al_{2}\mu\nu$,
$w_{1}=(\al_{1},\bt_{3}\nu)$;
\item[(4)]  $\al_{3}=\bt_{2}\mu$, $\al_{2}\mu=\bt_{1}\nu$,
$w_{1}=(\al_{1}\nu,\bt_{3})$;
\ei
One checks directly that in each of these four cases, $w_{2}=w_{1}$.  (In
case (4), one needs to consider separately the cases $\al_{2}=\bt_{1}\xi$
and $\bt_{1}=\al_{2}\xi$.)  So if $w_{1}\neq z$, then $w_{1}=w_{2}$.
Similarly, one shows that if $w_{2}\neq z$, then $w_{2}=w_{1}$.  The
associative law then follows.

Next we have to show that for each $t\in T$, $t^{*}$ is the only element
$t'$ for which $tt't=t, t'tt'=t'$.  If $t=z$, then this is trivial.  If
$t=(\al,\bt)$, then since $tt't=t=(\al,\bt)$, then $t'=(\al',\bt')$ for some
finite paths $\al',\bt'$, with $\bt',\al$ and $\bt,\al'$ comparable.
If $\bt'=\al\mu$ with $l(\mu)>0$,
then the second component of
$tt't$ would be longer than $\bt$, the second component of $t$,
and we contradict $tt't=t$.
Similarly, by considering $t'tt'=t'$ we obtain that $\al$ is not of the form
$\bt'\mu$ for any $\mu$ of positive length.
So $\bt'=\al$.  Similarly, $\al'=\bt$
and $t'=t^{*}$. Obviously $tt^{*}t=t,t^{*}tt^{*}=t^{*}$.
\end{proof}

\begin{proposition}   \label{prop:Phi}
The map $\Phi:S=S_{\calE}\to T$ where
\[    \Phi(v)=(v,v), \hspace{.1in}\Phi(e)=(e,r(e)), \hspace{.1in}
\Phi(e^{*})=(r(e),e), \hspace{.1in}\Phi(z)=z, \]
is a semigroup isomorphism.
\end{proposition}
\begin{proof}
Let $W=V\cup E\cup E^{*}\cup \{z\}$.
Then $\Phi$ extends in the obvious way to a
homomorphism, also denoted $\Phi: F_{W}\to T$.
The map $\Phi:F_{W}\to T$ is onto since every $t\in T$ is a
product of $\Phi(w)$'s: for example, if $\al=\al_{1}\ldots \al_{k}$ and
$\bt=\bt_{1}\ldots \bt_{l}$, and $t=(\al,\bt)\in T$, then
\[ (\al,\bt)= \Phi(\al_{1})\ldots \Phi(\al_{k})\Phi(\bt_{l}^{*})\ldots
\Phi(\bt_{1}^{*}).                \]
Let $\sim$ be the congruence on $F_{W}$ generated by the relations (i)-(iv)
above. To prove that $\Phi$ descends to a homomorphism $\tilde{\Phi}$
from $S=F_{W}/\sim$
onto $T$, we just have to check that the elements of $\Phi(W)$ satisfy
the relations in $T$ corresponding to (i)-(iv). This is straightforward.
(i) is trivially true. In (ii), for example, the equalities
\[  \Phi(s(e))\Phi(e)=\Phi(e)=\Phi(e)\Phi(r(e))   \]
are just
\[ (s(e),s(e))(e,r(e))=(e,r(e))=(e,r(e))(r(e),r(e)).                \]
In (iii), if $a=v\in V$ and
$r(a)\neq s(b)$, then either $b=w\in V$ with $w\neq v$, or $b=e$ with $v\neq
s(e)$ or $b=e^{*}$ with $v\neq r(e)$. Then $\Phi(a)\Phi(b)$ is either
$(v,v)(w,w)$ or $(v,v)(e,r(e))$ or $(v,v)(r(e),e)$ and these are all $z$. The
other possibilities for $a$ are dealt with similarly. For (iv),
$\Phi(e^{*})\Phi(f)=(r(e),e)(f,r(f))=z$ if $e\neq f$.

To show that $\tilde{\Phi}:S\to T$ is an
isomorphism, suppose that for some $s,s'\in F_{W}$
we have $\Phi(s)=\Phi(s')$. We have to show that $s\sim s'$.
Using the relations (i)-(iv), we can find $t,t'\in F_{W}$, where
$s\sim t$, $s'\sim t'$ and, either $t=t'=z$, in which case $s\sim s'$, or, in
an obvious notation, $t=\al\bt^{*}, t'=\al'(\bt')^{*}$ for paths
$\al,\bt,\al',\bt'$ with $r(\al)=r(\bt), r(\al')=r(\bt')$.  Then
\[ (\al,\bt)=\Phi(t)=\Phi(s)=\Phi(s')=\Phi(t')=(\al',\bt')  \]
and $\al=\al', \bt=\bt'$.  So $t=t'$ and $s\sim s'$.
\end{proof}

For the rest of the paper, we shall usually identify $S$ with $T$.

As a very simple example of the inverse semigroup $S$ of a graph $\calE$,
consider the case where $\calE$ has exactly one edge $e$ with initial and
terminal vertices $v_{1},v_{2}$ respectively.  (In this case, $\calE$ has a
sink.) The elements of $S$ are then:
$v_{1}=(v_{1},v_{1}), v_{2}=(v_{2},v_{2}), (e,v_{2}), (v_{2},e), (e,e)$ and the
zero $z$. Then the idempotent set $I$ of $S$ is
$\{v_{1}, v_{2}, (e,e),z\}$, and
$v_{1}v_{2}=z$, $v_{1}(e,e)=(e,e)$ and $v_{2}(e,e)=z$.
The involution on $S$
interchanges $(e,v_{2})$ and $(v_{2},e)$ and leaves every other element of
$S$ fixed. The other non-zero
products in $S$ are: $v_{1}(e,v_{2})=(e,v_{2})=(e,v_{2})v_{2},
(e,v_{2})(v_{2},e)=(e,e), (v_{2},e)(e,e)=(v_{2},e)$,
$(v_{2},e)(e,v_{2})=v_{2}$ and the equalities obtained from these by
applying the involution.

The next objective is to identify the universal groupoid of $S$ for
general $\calE$. The universal groupoid (\cite[Ch. 4]{Paterson})  $H$
of a countable inverse semigroup $S$ is constructed as follows. The unit
space $H^{0}$ of $H$ is just the set of non-zero semicharacters $\chi:I\to
\{0,1\}$, where $I$ is the commutative semigroup of idempotents of $S$. The
topology on $H^{0}$ is just the topology of pointwise convergence on $I$. There
is a dense subset $\ov{I}=\{\ov{\al}:\al\in I\}$ of $H^{0}$, where
$\ov{\al}(\bt)=1$ if
and only if $\al\bt=\al$. There is a natural right action of $S$
on $H^{0}$ given as follows.  First, an element
$x\in H^{0}$ is in the domain $D_{s}$ of $s\in S$ if $x(ss^{*})=1$. The
element $x.s\in H^{0}$ is given by: $x.s(\al)=x(s\al s^{*})$.
The map $x\to x.s$ is
a homeomorphism from $D_{s}$ onto $D_{s^{*}}$. The family of sets of the form
$D_{\al,\al_{1},\ldots \al_{n}}$, where
\begin{equation}
D_{\al,\al_{1},\ldots \al_{n}}=
D_{\al}\cap D_{\al_{1}}^{c}\cap \ldots \cap D_{\al_{r}}^{c} \label{eq:Uee}
\end{equation}
(with $c$ standing for ``complement'' and $\al,\al_{i}\in I, \al\geq \al_{i}, 1\leq i\leq r$) is a basis for the topology of $H^{0}$.

The universal groupoid $H$ is the quotient
       \[        \{(x,s): x\in D_{s}, s\in S\}/\sim      \]
where $(x,s)\sim (y,t)$ whenever $x=y$ and
there exists $\al\in I$ such that $x(\al)=1 \hspace{.1in}
(=x(ss^{*})=x(tt^{*}))$
and $\al s=\al t$.
The composable pairs are pairs of the form
$(\ov{(x,s)},\ov{(x.s,t)})$ ($x\in D_{s}$, $x.s\in D_{t}$, $s,t\in S$), while
the product and involution on $H$ are respectively given by
$(\ov{(x,s)},\ov{(x.s,t)})\ra \ov{(x,st)}$ and $\ov{(x,s)}\ra
\ov{(x.s,s^{*})}$. Then $H$ is an r-discrete groupoid and the map $\Psi$,
where $\Psi(s)=\{\ov{(x,s)}:x\in D_{s}\}=A_{s}$ is an inverse semigroup
isomorphism from $S$ into the ample semigroup $H^{a}$ of $H$. (Recall that for
any r-discrete groupoid $G$, the {\em ample semigroup} $G^{a}$ is the inverse
semigroup of compact, open, Hausdorff subsets of $G$.)  Note that
$\Psi(\al)=D_{\al}$ for all $\al\in I$. A basis for the
topology of $H$ is given by sets of the form $D_{\al,\al_{1},\ldots
,\al_{n}}\Psi(s)$ ($\al,\al_{i}\in I, \al_{i}\leq \al, s\in S$).

We now specify the universal groupoid $H=H_{\calE}$ of
$S=S_{\calE}$. The discussion
parallels that for the Cuntz inverse semigroups $S_{n}$ ($1\leq n\leq\infty$)
given in \cite[pp.182f.]{Paterson}.  We first determine the unit space
$H^{0}$ of $H$.  For $\al\in I$, we will write $\al$ in place of $\ov{\al}$.
The context will make clear if $\al$ is being regarded as an element of
$I$ or as a semicharacter. Recall that $Y$ and $Z$ are respectively the
sets of finite and infinite paths in $\calE$.

\begin{proposition}   \label{prop:unit}
The unit space $H^{0}$ of $H$ is the disjoint union:
\begin{equation}
	H^{0}=Y\cup Z\cup\{z\}.    \label{eq:Ghat0}
\end{equation}
Then the topology for $H^{0}$ is specified by:
\bi
\item[(i)] the singleton set $\{z\}$ is both open and closed in $H^{0}$;
\item[(ii)] a neighborhood basis at $u\in Z$ is given by sets of the form
$D_{\al} = \{\al \ga\in Y\cup Z:\ga\in Y\cup Z\}$ ($\al\in Y$) where $\al$ is
an initial segment of $u$;
\item[(iii)] if $y\in Y_{f}$, then $\{y\}$ is an open set;
\item[(iv)] if $y\in Y_{\infty}$, then a neighborhood basis at $y$
is given by sets of the form $D_{y,ye_{1}, \dots ,ye_{n}}$ where
$e_{1},\ldots e_{n}$ are edges in $\calE$ starting at $r(y)$.
\ei
\end{proposition}
\begin{proof}
Let $\chi\in H^{0}$. If $\chi(z)=1$, then since $\al z=z$ for all $\al\in I$,
we have $\chi(\al)=1$, the constant non-zero semicharacter
on $E$. Clearly $\chi=z$ $(=\ov{z})$.
The set $\{z\}$ is open in $H^{0}$ since it is just $\{\chi'\in
H^{0}: \chi'(z)=1\}$. This gives (i).

Suppose then that $\chi\neq z$.  Then $\chi(z)=0$.  Recall that $I$ is
identified with the set of finite paths $\al$
together with the zero $z$, and a
product of two of these paths is $\neq z$ only when one path is an initial
segment of the other, and in that case, the product is just the longer path.
The set of paths $\al$ for which $\chi(\al)=1$ is a subsemigroup of $I$ and
so each such path is an initial segment of any other or vice versa.  There
thus exists a path $\ga\in Y\cup Z$ such that $\chi(\al)=1$ if and only
if $\ga$ is of the form $\al\mu$. Identifying $\chi$ with $\ga$ gives
(\ref{eq:Ghat0}). Note that if $\ga\in Y$, then $\chi=\ga$ ($=\ov{\ga}$).

By the above, a basis for the topology on $Y\cup Z$ is given by sets of the
form $D_{\al,\al_{1}, \cdots ,\al_{n}}$. Suppose that
$u\in Z$ and that $u\in D_{\al,\al_{1},
\cdots ,\al_{n}}$. Then $u$ starts with $\al$ but not with any $\al_{i}$.  So
for some $n$, $u_{1}\cdots u_{n}\in D_{\al,\al_{1},
\cdots ,\al_{n}}$, and so $D_{u_{1}u_{2}\cdots u_{n}}= \{u_{1}\cdots
u_{n}\ga\in Y\cup Z: \ga\in Y\cup Z\} \subset D_{\al,\al_{1}, \cdots
,\al_{n}}$. So the sets of the form $D_{\al}$ ($\al=u_{1}\cdots u_{n}$) form a
neighborhood base at $u$ in $Y\cup Z$. This gives (ii).

For (iii), if $y\in Y_{f}$ and
$e_{1},\ldots ,e_{n}$ are all of the edges starting at $r(y)$, then
\[ \{y\}=D_{y}\cap D_{ye_{1}}^{c}\cap\cdots \cap D_{ye_{n}}^{c}=
D_{y,ye_{1},\ldots
,ye_{n}}  \]
is open in $Y\cup Z$.

Lastly, suppose that $y\in Y_{\infty}$ and that $y$ is in
some $D_{\al,\al_{1},\ldots ,\al_{n}}$. Then $y$ begins with $\al$ and the
$\al_{i}$'s all begin with $y$ and are longer than $y$. Write
$\al_{i}=ye_{i}\ldots$. Then $y\in D_{y,ye_{1}, \dots ,ye_{n}}\subset
D_{\al,\al_{1},\ldots ,\al_{n}}$ and (iv) follows.
\end{proof}

It is instructive to specify the convergent sequences in $Y\cup Z$. If $u\in
Z$, then $x^{n}\to u$ if and only if $l(x^{n})\to \infty$, and for each $i$,
$x^{n}_{i}=u_{i}$ eventually. Next, let us say that a sequence $\{e^{n}\}$ in
$E$ is {\em wandering} if, for any $e\in E$, $e^{n}\neq e$ eventually. Let
$\al=\al_{1}\al_{2}\ldots \al_{N} \in Y$. Suppose that $x^{n}\in Y\cup Z$
and that $x^{n}\neq \al$ for any $n$. Then $x^{n}\to \al$ if and only if,
eventually, $l(x^{n})> l(\al)$, $x^{n}_{i}=\al_{i}$ for $1\leq i\leq N$, and
$\{x^{n}_{N+1}\}$ is wandering. Note that in this case, $\al\in Y_{\infty}$. It
follows that each $y\in Y_{\infty}$ is a limit point of $Z$, and so we have the
following proposition.

\begin{proposition}          \label{prop:dense}
The set $Z$ of infinite paths is dense in $Y_{\infty}\cup Z$. If $\al\in Y$,
then $\{\al\}$ is open in $H^{0}$ if and only if $\al\in Y_{f}$.
\end{proposition}

\begin{theorem}    \label{th:univ}
The universal groupoid $H$ for $S$ can be identified with the union
of $\{z\}$ and the set of
all triples of the form $(\al\ga,l(\al)-l(\bt),\bt\ga)$ where $\al,\bt\in
Y$, $\ga\in Y\cup Z$, and $\al\ga, \bt\ga\in Y\cup Z$. Multiplication on
$H$ is given by: $(x,m,y)(y,n,z)=(x,m+n,z)$ and inversion by
$(x,m,y)^{-1}=(y,-m,x)$. The canonical map
$\Psi:S\to H^{a}$ sends $z$ to $\{z\}$ and any element $(\al,\bt)$ to the
(compact open) set $A_{\al,\bt}$ of all
triples of the form $(\al\ga,l(\al)-l(\bt),\bt\ga)$. The locally compact
groupoid $H$ is Hausdorff.
\end{theorem}
\begin{proof} The proof is close to
that for the Cuntz semigroup $S_{n}$ in \cite[pp.182-186]{Paterson}. We can
restrict consideration to $H\sim \{z\}$. Let $s=(\al,\bt)$, $x\in D_{s}$. This
means that $x(ss^{*})=x(\al)=1$. If $\al'\in I$ is such that
$x(\al')=1, \al'\leq
ss^{*}$, then $\al'=\al\de\in Y$, and $\al's=(\al\de,\bt\de)$.
Then $x=\al\de\ga$ for some $\ga$ and
$(\al\de\ga,l(\al\de)-l(\bt\de),\bt\de\ga)=
(\al\de\ga,l(\al)-l(\bt),\bt\de\ga)$. So the map $(x,s)\to
(\al\de\ga,l(\al)-l(\bt),\bt\de\ga)$ is constant on the equivalence class
of $(x,s)$ in $H$, and so defines a map on $H$. The argument reverses to give
that this map is a bijection. We now prove that $H$ is Hausdorff, leaving the
remaining verifications of the theorem to the reader.

Let $a=(\al, l(\al)-l(\bt),\bt)$, $b=(\al',l(\al')-l(\bt'),\bt')$ belong to
$H$ with $a\neq b$. If $A_{\al,\bt}\cap A_{\al',\bt'}=\emptyset$, then
then we can separate $a$ and $b$ using $A_{\al,\bt}, A_{\al',\bt'}$. Suppose
then that $A_{\al,\bt}\cap A_{\al',\bt'}\neq \emptyset$.
Then there exist $\ga,\ga'$ such that
\[  (\al\ga,l(\al)-l(\bt),\bt\ga)=(\al'\ga',l(\al')-l(\bt'),\bt'\ga').   \]
Then $l(\al)-l(\bt)=l(\al')-l(\bt')$ and $\al,\al'$ and $\bt,\bt'$ are
comparable. We can suppose that for some $\mu$, $\al'=\al\mu$. Then
$\bt'=\bt\nu$ where $l(\mu)=l(\nu)$, and $\nu=\mu$ since both are initial
segments of $\ga$. Since $a\neq b$, we have $l(\mu)>0$. Then $A_{\al,\bt}\cap
A_{\al',\bt'}^{c}$, $A_{\al',\bt'}$ separate $a$ and $b$.
So $H$ is Hausdorff.
\end{proof}

We next set $X=Y_{\infty}\cup Z\subset H^{0}$. Then
(Proposition~\ref{prop:unit}) $X$ is a closed invariant subset of $H^{0}$
and so is a locally compact Hausdorff space in its own right.
Let $G_{\calE}=G$ be the
reduction of $H$ to $X$. We will call $G$ the {\em path groupoid} of $\calE$.
(The groupoid $G$ coincides with the usual path groupoid of a
graph in the locally finite case (e.g. \cite{KPRR}).)  Then $G$ is a
closed subgroupoid of $H$, and is an r-discrete groupoid in its own right
with counting measures giving a left Haar system. For $\al,\bt\in Y$, we
define $A'(\al,\bt)=A(\al,\bt)\cap G$ and $D_{\al}'=D_{\al}\cap X$.

We now relate the representations of $G$ to those of $H$.

\begin{proposition}          \label{prop:rep}
Each representation $\si$ of $C_{c}(G)$ determines a representation $\si'$ of
$C_{c}(H)$ by: $\si'(f)=\si(f\mid G)$. A representation $\Pi$ of $C_{c}(H)$
is of the form $\si'$ for some representation $\si$ of $C_{c}(G)$
if and only if
\begin{equation}
		\Pi(\chi_{\{z\}})=0=\Pi(\chi_{\{\al\}})   \label{eq:zal}
\end{equation}
for all $\al\in Y_{f}$.
\end{proposition}
\begin{proof}
If $\si$ is a representation of $C_{c}(G)$ on a Hilbert space $\mfH$
then $\si':C_{c}(G)\to B(\mfH)$ is a $^{*}$-homomorphism. Since the
homomorphism $f\to
f\mid G$ from $C_{c}(H)$ into $C_{c}(G)$ is continuous for the $I$-norm, it
immediately follows that $\si'$ is $I$-norm continuous and so is a
representation of $C_{c}(H)$. It is obvious from the definition that $\si'$
satisfies (\ref{eq:zal}) (with $\si'$ in place of $\Pi$).

Conversely, suppose that $\Pi:C_{c}(H)\to B(\mfH)$ is a representation
satisfying (\ref{eq:zal}).  By the proof of
\cite[Theorem 3.1.1]{Paterson}, the
representation $\Pi$ extends to a continuous
representation $\tilde{\Pi}$ of
$(B_{c}(H),\norm{.}_{I})$ where $B_{c}(H)$ is the convolution algebra of
bounded Borel functions on $H$ with compact supports.  Let $\si$ be the
restriction of $\tilde{\Pi}$ to $C_{c}(G)\subset B_{c}(H)$.  Then $\si$ is a
representation ({\em a priori} possibly degenerate) of $C_{c}(G)$ on $\mfH$.
Let $f\in C_{c}(H)$ and $U$ be an open subset of $H$ with compact closure and
which contains the support of $f$.  We have to show that
$\Pi(f)=\tilde{\Pi}(f\mid G)$ $(=\si(f\mid G)).$
Let $\xi\in \mfH$.  As in
\cite[Theorem 3.1.1]{Paterson}, there exists a positive regular Borel
measure $\mu$ on $H$ such that for all $g\in B_{c}(H)$ with
$supp(g)\subset U$, we have
\begin{equation}
\langle\tilde{\Pi}(g)\xi,\xi\rangle=\int_{U}g\,d\mu.  \label{eq:star}
\end{equation}
Now $H\setminus G$ is the reduction of $H$ to $Y_{f}\cup \{z\}$ and is
countable.  For $a\in H\setminus G$, we have $a=ar(a)$ with $r(a)\in Y_{f}\cup
\{z\}$.  So $\Pi(\chi_{\{a\}})=\Pi(\chi_{\{a\}})\Pi(\chi_{\{r(a)\}})=0$ by
(\ref{eq:zal}).  From (\ref{eq:star}), $\mu(\{a\})=0$.  Since $H\setminus G$
is countable, it follows that $\mu(H\setminus G)=0$, and
for $f\in C_{c}(G)$,
\[\langle\tilde{\Pi}(f)\xi,\xi\rangle
=\int_{U} f\mid G\,d\mu
=\langle\tilde{\Pi}(f\mid G)\xi,\xi\rangle,       \]
and $\Pi(f)=\si(f\mid G)$.
\end{proof}

The requirement (\ref{eq:zal})
of the preceding proposition can be replaced by
the corresponding requirement in which $\al\in Y_{f}$ is replaced by $v\in
V_{f}$:
\begin{equation}
 \Pi(\chi_{\{z\}})=0=\Pi(\chi_{A(v,v)}) - \sum_{s(e)=v}\Pi(\chi_{A(e,e)}).
		  \label{eq:pizal2}
\end{equation}
To show this, $\Pi(\chi_{A(v,v)}) -
\sum_{s(e)=v}\Pi(\chi_{A(e,e)})=\Pi(\chi_{\{v\}})$ and
one way follows trivially since $v\in Y_{f}$. For the other way, we just use
the fact that $\chi_{\{\al\}}=
\chi_{\{\al\}}*\chi_{\{r(\al)\}}$ and the fact that $\Pi$ is a
homomorphism.

We now give the definition of a {\em representation} of the graph $\calE$. The
definition is taken from \cite{FLR}. Note that in the $\mathcal{O}_{\infty}$
case, condition (iii) below assumes the form (\cite[p.174]{CuntzSimple})
$\sum_{i=1}^{r}S_{i}S_{i}^{*}\leq 1$ for all $r\geq 1$. Further, $1$ is just
$P_{v}$ where $v$ is the sole vertex of the directed graph which is the
infinite bouquet of circles.

\begin{definition}    \label{def:rep}
A {\em representation} of $\mathcal{E}$ on a Hilbert space $\mfH$ is given by
a family $\{P_{v}: v\in V\}$
of mutually orthogonal projections,
and a family $\{S_{e}:e\in E\}$
of partial isometries such that the family of projections
$S_{e}S_{e}^{*}$ are mutually orthogonal and:
\be
\item[(i)] $S_{e}^{*}S_{e}=P_{r(e)}$ for all $e\in E$;
\item[(ii)] if $v\in V_{f}$, then
\begin{equation}
P_{v}=\sum_{s(e)=v}S_{e}S_{e}^{*};         \label{eq:pvsef}
\end{equation}
\item[(iii)] $P_{s(e)}S_{e}S_{e}^{*}=S_{e}S_{e}^{*}$ for all $e\in E$.
\ee
\end{definition}

By multiplying both sides of (i) on the right by $S_{e}^{*}$ we get
\begin{equation}
S_{e}^{*}=P_{r(e)}S_{e}^{*}.  \label{eq:se1}
\end{equation}
Similarly, using (iii), we get that
\begin{equation}
P_{s(e)}S_{e}=S_{e}.  \label{eq:se2}
\end{equation}
In the following, representations of graphs and inverse
semigroups are assumed to be non-degenerate in the obvious sense.
(Representations of locally compact groupoids are ``non-degenerate'' almost by
definition.)

\begin{theorem} \label{th:rep}
There are natural one-to-one correspondences between:
(a) the class $\mfR_{\calE}$ of representations of $\calE$;
(b) the class $\mfR_{S}$ of representations of $S$
that vanish on $z$ and on elements of the form $v-\sum_{s(e)=v}(e,e)$ for all
$v\in V_{f}$;
(c) the class $\mfR_{G}$ of representations of the locally compact
groupoid $G$.
\end{theorem}
\begin{proof}
Let $Q=\{\{P_{v}\}, \{S_{e}\}\}\in \mfR_{\calE}$ and be realized on a Hilbert
space $\mfH$. For $\al=e_{1}e_{2}\ldots e_{n}\in Y$, define
$S_{\al}=S_{e_{1}}S_{e_{2}}\ldots S_{e_{n}}$. (i) above generalizes to:
\begin{equation}
S_{\al}^{*}S_{\al}=P_{r(\al)}.  \label{eq:salpha}
\end{equation}
For example, if $\al=e_{1}e_{2}$, then using (i), (iii) and (\ref{eq:se2}),
we have
$S_{\al}^{*}S_{\al}=S_{e_{2}}^{*}[S_{e_{1}}^{*}S_{e_{1}}]S_{e_{2}}
=S_{e_{2}}^{*}P_{s(e_{2})}S_{e_{2}}
=(S_{e_{2}}^{*}S_{e_{2}})=P_{r(e_{2})}=P_{r(\al)}.$
It immediately follows from (\ref{eq:se2}) and (\ref{eq:se1}) that:
\begin{equation}
P_{s(\al)}S_{\al}=S_{\al}, P_{r(\al)}S_{\al}^{*}=S_{\al}^{*}. \label{eq:psalpha}
\end{equation}

Define a $^{*}$-map $\Pi:S\to B(\mfH)$ by:
\[  \Pi(v)=P_{v},\hspace{.1in} \Pi(\al,\bt)=S_{\al}S_{\bt}^{*},
\hspace{.1in}\Pi(z)=0. \]
Then $\Pi$ is a $^{*}$-homomorphism. The proof is simple.  One
checks that $\Pi(st)=\Pi(s)\Pi(t)$ for the different kinds of product using
(\ref{eq:salpha}) and (\ref{eq:psalpha}).  For example, to prove that
$\Pi(\al,\bt\mu)=\Pi(\al,\al'\mu)\Pi(\al',\bt)$, one argues: \\
$\Pi(\al,\al'\mu)\Pi(\al',\bt)
=S_{\al}S_{\mu}^{*}S_{\al'}^{*}S_{\al'}S_{\bt}^{*}
=S_{\al}S_{\mu}^{*}P_{r(\al')}S_{\bt}^{*}
=S_{\al}S_{\mu}^{*}S_{\bt}^{*}=\\
\Pi(\al,\bt\mu).$ Since $\Pi(z)=0$, $S_{e}S_{e}^{*}=\Pi((e,e))$ and $\Pi$
vanishes on elements of the form $(v-\sum_{s(e)=v}(e,e))$ for $v\in V_{f}$, it
follows from (\ref{eq:pvsef}) that $\Pi\in \mfR_{S}$.

Conversely, any $\Pi'\in \mfR_{S}$ determines an element of $\mfR_{\calE}$ by
taking $P_{v}=\Pi'(v,v), S_{e}=\Pi'(e,r(e))$. Since
$S_{e}S_{e}^{*}=\Pi'((e,r(e))(r(e),e))= \Pi'(e,e)$, (ii) of
Definition~\ref{def:rep} follows. (i) and (iii) follow since $\Pi'$ is a
$^{*}$-homomorphism on $S$. This establishes the correspondence between the
representations of (a) and (b).

The correspondence between the representations in (b) and (c) is implemented
by using the natural bijection between representations of $S$ and
$(H,\Psi)$ (the universal groupoid of $S$) established in
\cite[Ch.4]{Paterson}. We now recall what this bijection is.

Let $\rho$ be a representation of $S$ on a Hilbert space $\mfH$.  Then there
is a representation $\Pi$ of $C_{c}(H)$ determined by:
\[  \Pi(\chi_{\Psi(s)})=\rho(s) \]
for all $s\in S$ where $\Psi(s)$ is given in Theorem~\ref{th:univ}.
Further, the map $\rho\to \Pi$ is one-to-one and onto.
It remains then to show that the
representations $\rho$ of $S$ that satisfy the conditions of (b)
correspond to the representations $\Phi$ of $C_{c}(G)$.
This follows by Proposition~\ref{prop:rep} (and the remark following that
proposition), and so we have established the
correspondence between $\mfR_{S}$ and $\mfR_{G}$.
\end{proof}

We now translate (a), (b) and (c) of the preceding theorem into $\cs$ terms.
The \cs\ determined by the groupoid representations of (c) is just the
groupoid \cs\ $C^{*}(G)$. The \cs\ associated with (b) is obtained on the
semigroup algebra $\ell^{1}(S)$ as follows. Every representation $\Pi$ of $S$
gives a bounded representation $\Pi$ of $\ell^{1}(S)$ in the natural way, and
$C^{*}(S)$ is just the enveloping $\cs$ of $\ell^{1}(S)$ obtained by taking the
biggest $C^{*}$-norm coming from such $\Pi$'s. The $\cs$ that we want here,
which we will denote by $C_{0}^{*}(S)$, is obtained in the same way only
using $\Pi$'s for which $\Pi(z)=0=\Pi(v-\sum_{s(e)=v}(e,e))$ for all $v\in
V_{f}$. It is easy to show that $C_{0}^{*}(S)$ is the enveloping \cs\ of
$\ell^{1}(S)/I$, where $I$ is the closed ideal of $\ell^{1}(S)$ generated by
elements of the form $z$ and $v-\sum_{s(e)=v}(e,e)$ ($v\in V_{f}$).

We take $C^{*}(\calE)$ to be the universal $\cs$ for the representations of
$\calE$. By the above, we then obtain the
following corollary to Theorem~\ref{th:rep}.

\begin{corollary} \label{cor:rep}
$C^{*}(\calE)\cong C_{0}^{*}(S)\cong C^{*}(G).$
\end{corollary}

In the case of $\mathcal{O}_{\infty}$, the graph $\calE$ is a bouquet of
circles $\{e_{n}:n\geq 1\}$ (one vertex $v$). The inverse semigroup $S_{\calE}$
in this case is the set of pairs $(\al,\bt)$ where $\al,\bt$ are finite words
in $\Pos$ together with the identity $v$ and zero $z$. Then $G$ is the groupoid
$O_{\infty}$ of \cite[p.184f.]{Paterson}.

Now let $\calE$ be a general directed graph and $S=S_{\calE}$ as before.
For the purposes of the next section, we need to use a related inverse
semigroup $S^{n}$ with associated groupoid $G^{n}$.  Here, $n\in \Pos$,
$Y^{n}=\{y\in Y: 0\leq l(y)\leq n\}$ and
$S^{n}$ is the inverse subsemigroup of $S$ given by:
  \[S^{n}=\{(\al\mu,\bt\mu)\in S: \al,\bt\in Y^{n}, l(\al)=l(\bt)\}\cup \{z\}.\]
It is easy to check that $S^{n}$ is an inverse subsemigroup of $S$ and that
$I(S^{n})=I(S)$.  We now describe the universal groupoid $H^{n}$ of $S^{n}$
and the related groupoid $G^{n}$. The proofs of the assertions below are very
similar to the corresponding proofs for $S$, and are left to the reader.

The universal groupoid $H^{n}$ for $S^{n}$ can be identified with the union of
$\{z\}$ and the set of all pairs of the form $(\al\ga,\bt\ga)$ where
$\al,\bt\in Y^{n}$, $l(\al)=l(\bt)$,
$\ga\in Y\cup Z$, and $\al\ga, \bt\ga\in Y\cup Z$.
Identify $(\al\ga,\bt\ga)\in S^{n}$ with $(\al\ga,0,\bt\ga)\in H$. Then $H^{n}$
is a subgroupoid of $H$ with the same unit space, and the multiplication and
inversion for $H^{n}$ are just those that $H^{n}$ inherits as a subgroupoid of
$H$. The canonical map $\Psi^{n}:S^{n}\to (H^{n})^{a}$ is just the restriction
of the map $\Psi$ of Theorem~\ref{th:univ} to $S^{n}$. The groupoid $G^{n}$ is
defined to be the reduction of $H^{n}$ to $X$.

The principal transitive equivalence relation on a (discrete) set with $k$
elements ($1\leq k\leq \infty$) is denoted by $T_{k}$. (In less elaborate
language, this is the equivalence relation for which any element of the set is
equivalent to any other.) For the following definitions, see \cite[p.123]{rg}.
An r-discrete groupoid is called {\em elementary} if it is
isomorphic to the
disjoint union of a sequence of groupoid products of the form $T_{k}\x Q$ for
some $k$ and some locally compact space $Q$, and is called {\em AF} if its unit
space is totally disconnected and it is the inductive limit of a sequence of
elementary groupoids.

\begin{proposition}   \label{prop:AF}
The locally compact groupoid $G^{n}$ is an AF-groupoid.
\end{proposition}
\begin{proof}
Obviously, the unit space $G^{0}$ of $G^{n}$ is totally disconnected.
Enumerate $E$ as $e_{1},e_{2},\ldots$, and let $N\in \Pos$. Let
$F=\{e_{1},\ldots ,e_{N}\}$. For $0\leq r\leq n$, let $A_{r}$ be the set of all
finite paths of the form $e_{p_{1}}\ldots e_{p_{r}}$ where each $e_{p_{i}}\in
F$. We take $A_{0}=s(F)$. Let $B_{r}$ be the set of
pairs $(\al\ga,\bt\ga)\in G^{n}$ where $\al,\bt\in A_{r}$ and $\ga$ is
arbitrary. Clearly, $B_{r}$ is a compact, open subgroupoid of $G^{n}$ (being
a finite union of $\Psi^{n}(\al,\bt)$'s). For each $v\in r(F)$, let
	     \[           X^{v}=\{\ga\in X: s(\ga)=v\}.           \]
Then $X^{v}=A(v,v)$ is a compact, open subset of $X$.  Note that if $\al,\bt\in
A_{r}$ with $r(\al)=v=r(\bt)$, then $(\al\ga,\bt\ga)\in B_{r}$ if and only if
$\ga\in X^{v}$.  Let
\[   L^{N}=\cup_{r=0}^{n} B_{r}.      \]
Then $L^{N}$ is also a compact, open subgroupoid of $G^{n}$.  We will show
that $L^{N}$ is elementary.

To this end we ``disjointify'' the $B_{r}$'s.  Let $C_{r}$ be the unit space
$B_{r}^{0}$ of $B_{r}$.  The elements of $C_{r}$ are of the form $\al\ga$
for $\al\in A_{r}, \ga\in X^{r(\al)}$.  Clearly
$C_{n}\subset C_{n-1}\subset \cdots \subset C_{0}=\{x\in X: s(x)\in A_{0}\}$.
Let $C_{n}'=C_{n}$ and, for $0\leq r< n$, let $C_{r}'=C_{r}\setminus C_{r+1}$.
Let
\[  B_{r}'=\{(\al\ga,\bt\ga)\in B_{r}: \al,\bt\in A_{r},
\al\ga, \bt\ga\in C_{r}'\}.    \]
Then $L^{N}$ is the disjoint union of the compact, open subgroupoids $B_{r}'$
of $G^{n}$.  The fact that $L^{N}$ is elementary will follow once we have
shown that each $B_{r}'$ is elementary.  To this end, fix $r$, and
for each $v\in r(F)$, let
\[ W_{v}=\{(\al\ga,\bt\ga)\in B_{r}': \al,\bt\in A_{r},
v=r(\al)=r(\bt)=s(\ga)\}.\]
Then $B_{r}'$
is the disjoint union of the compact, open subgroupoids $W_{v}$, and for each
$v$, $W_{v}$ is isomorphic to $T_{k_{v}}\x (X^{v})'$, where $(X^{v})'$ is an
open subset of $X^{v}$, and
$k_{v}$ is the
number of $\al$'s in $A_{r}$ with $r(\al)=v$ and
such that $\al\ga\in C_{r}'$ for some $\ga$. So $B_{r}'$, and hence $L^{N}$, is
elementary.

We now claim that $G^{n}$ is the inductive limit of the $L^{N}$'s and hence
is an AF-groupoid.  To this end (cf. \cite[pp.122-3]{rg}), $G^{n}$ is a union
of the increasing sequence of elementary, open subgroupoids $L^{N}$.  It is
easily checked that a set $V\subset G^{n}$ is open in the inductive limit
topology if and only if it is open in $G^{n}$, so that the inductive limit
topology on $G^{n}$ coincides with the given topology on $G^{n}$.  Further,
the counting measure left Haar system on $G^{n}$ is compatible with the
counting measure left Haar systems on the $L^{N}$'s.  So $G^{n}$ is the
inductive limit of the $L^{N}$'s.
\end{proof}

At the $S^{n}$ level (for $\mathcal{O}_{\infty}$), Cuntz
observes (\cite[1.5]{CuntzSimple}) that $\mathcal{F}^{\infty}$ is an AF $\cs$.


\section{Graph groupoids are amenable}
In \cite{CuntzSimple}, Cuntz showed that $\mathcal{O}_{n}$ ($n\geq 2$) is
isomorphic to the crossed product of an AF-algebra by an automorphism, cut down
by a projection. It follows (\cite{CuntzSimple,rosen}) that
$\mathcal{O}_{n}$ is therefore nuclear (=amenable). Cuntz used a separate
argument to deal with the case $n=\infty$. Renault (\cite[p.138f.]{rg}) showed
that the Cuntz groupoids $O_{n}$ ($n<\infty$) are amenable as locally compact
groupoids. (The amenability of $\mathcal{O}_{n}$ then follows.) The groupoid
case when $n=\infty$ remained unclear.

Kumjian, Pask, Raeburn and Renault (\cite{KPRR}) adapted Renault's original
argument to establish the amenability of the groupoid $G=G_{\calE}$ when
$\calE$ is locally finite. In this section, we further adapt the argument of
\cite{KPRR} to show that $G$ is amenable in complete generality. In particular,
this covers the groupoid case for $O_{\infty}$.

As in \cite{KPRR}, we assume initially that every vertex receives an edge. For
each $v\in V$, fix an edge $e(v)\in E$ for which $r(e(v))=v$. Next we
consider the following space $\mfX$ of two-sided paths.
The elements of $\mfX$ are sequences $x=\{x_{i}\}$ of edges, where
$-\infty<i<k\leq\infty$ with $k$ depending on $x$ and
$r(x_{i})=s(x_{i+1})$ whenever that makes sense. So the sequence $x$ is
infinite to the left, and may or may not be infinite to the right. If the
sequence $x$ is finite to the right, then $r(x)$ makes sense in the obvious
way, and we require that $r(x)\in V_{\infty}$. Define a function $h:\mfX\to
\Z\cup \{\infty\}$ as follows: $h(x)=\infty$ if $x$ is infinite to the
right, and otherwise, is the largest $i$ for which $x_{i}$ is defined.

For each $n\in \Z$, we now define a subset $P_{n}$ of $\mfX$ by specifying
that $x\in P_{n}$ if and only if $h(x)\geq (n-1)$ and $x_{i}=e(s(x_{i+1}))$
for all $i<n$. If $h(x)=n-1$, then we require $x_{n-1}$ to be some $e(v)$.
The set of $x$'s in $P_{n}$ for which $h(x)=n-1$ can be identified with
$V_{\infty}$ by sending $x$ to $r(x)$. Note that $P_{n}\subset P_{n-1}$.

Let $n\geq 1$. There is a natural bijection $f_{-n}:P_{-n}\to X=G^{0}$ defined
by translating to the right by $n+1$ and chopping off what is to the left of
$1$. Precisely, $f_{-n}(x)_{m}=x_{-n+m-1}$ if $h(x)>(-n-1)$ and
$f_{-n}(x)=r(x)$
when $h(x)=-n-1$. We give $P_{-n}$ the topology for which $f_{n}$ is a
homeomorphism.

Note that $P_{-n}$ is an open subset of $P_{-n-1}$.  Indeed $f_{-n-1}(P_{-n})$
is the (open) set of elements of $X$ of the form $e(s(x))x$ for some
$x\in X$.

Let $R_{n}$ be the equivalence relation on $P_{-n}$ defined by: $(x,y)\in
R_{n}$ if and only if $h(x)=h(y)$, and if $h(x)=-n-1$ then $r(x)=r(y)$, while
if $h(x)>n$, then $x_{i}=y_{i}$ for $i>n$. Note that when $h(x)=-n-1$, then
$xR_{n} y$ if and only if $x=y$. Note also that if $h(x)>n$, then since
$x_{n+1}=y_{n+1}$, we have $r(x_{n})=s(x_{n+1})=s(y_{n+1})=r(y_{n})$.
The unit space of $R_{n}$ is $P_{-n}$.  Of course, we give $R_{n}$ the
product topology that it inherits as a subset of $P_{-n}\x P_{-n}$.

\begin{proposition}   \label{prop:Rn}
The equivalence relation groupoid $R_{n}$  is an AF groupoid.
\end{proposition}
\begin{proof}
The map $(x,y)\to (f_{n}(x),f_{n}(y))$ is an isomorphism of locally compact
groupoids from $R_{n}$ onto $G^{2n+1}$ and therefore is an AF groupoid by
Proposition~\ref{prop:AF}.
\end{proof}


Clearly, $R_{n}$ is an open subset of $R_{n+1}$, and the topology on $R_{n}$ is
the relative topology that it inherits as a subset of $R_{n+1}$. Also, the
canonical left Haar system $\{\la_{n+1}(u)\}$ (counting measure on the
$(R_{n+1})^{u}$'s) gives by restriction the left Haar system $\{\la_{n}(u)\}$ on
$R_{n}$.

It follows as in \cite[p.122]{rg} that $R=\cup_{n\geq 1}R_{n}$ is an r-discrete
groupoid in the inductive limit topology. (As pointed out in \cite{KPRR}, this
notion of inductive limit is more general than the r-discrete case
considered by Renault inasmuch as the unit spaces of the $R_{n}$'s are
increasing rather than being fixed.) It also follows that $R$ is amenable
since each $R_{n}$ is AF and therefore amenable (\cite[p.123]{rg}).

The argument for the amenability of $G$ now proceeds almost exactly as in
\cite{KPRR} (to which the reader is referred for further details). We will be
content with a sketch of the argument in the present context to indicate the
minor changes needed.

Let
$P=\cup_{1}^{\infty} P_{-n}$ be the unit space of $R$. Then $P$ is the
inductive limit of the $P_{-n}$'s. The left shift map $h$ on $P$, where
$h(\{x_{i}\})_{j} =x_{j+1}$, is a homeomorphism on $P$ since it maps $P_{-n}$
homeomorphically onto $P_{-n-1}$, and it induces a homeomorphism $\si$ on $R$,
where $\si(u,v)=(h(u),h(v))$, which is also a groupoid homomorphism. One then
considers the
semi-direct product $R\x_{\si} \Z$.  This groupoid consists (\cite{KPRR})
of triples $(u,v,k)\in R\x \Z$ with $r(u,v,k)=u, s(u,v,k)=h^{-k}v$, inverse by
$(u,v,k)^{-1}=(\si^{-k}(v,u),-k)$ and product by:
\[  (u,h^{k}(p),k)(p,q,l)=(u,h^{k}(q),k+l).      \]
Then $R\x_{\si} \Z$ is amenable (\cite[p.96]{rg}) and hence (\cite[p.92]{rg})
so is its reduction $G'$ to the closed subset $P_{1}$ of $P$.

We now claim that $G'$ is isomorphic to the graph groupoid $G$ (so that $G$ is
also amenable). Firstly, the unit space $P_{1}$ of $G'$ can be identified with
$X$ by sending $\ldots  e(s(z_{1}))z_{1}z_{2}\ldots$ to $z=z_{1}z_{2}\ldots$.
The natural isomorphism $F:G\to G'$ is given by
(cf. \cite[Proposition 4.4]{KPRR}):
\[F(xz,k,yz)= (xz,h^{-k}(yz),-k).   \]
Here, $x=x_{1}\cdots x_{r}, y=y_{1}\cdots y_{r-k}$ and $xz,yz$  on the
right-hand side are abbreviations for $\cdots e(s(x_{1}))xz,
\cdots e(s(y_{1}))yz\in P_{1}$.

Checking that $F$ is a groupoid isomorphism is straighforward.  That $F$ is a
homeomorphism follows effectively from the way that we defined the topology
of $P_{-n}$ in terms of $X$.  To check that $F$ is a homomorphism, one
argues ($k=l(x)-l(y), l=l(y_{1})-l(a), yz=yy'z'=y_{1}z'$):
\beqns
F(xz,k,yz)F(yz,l,az') & = & (xz,h^{-k}(yz),-k)(yz,h^{-l}(az'),-l)  \\
& = & (xz,h^{-k}[h^{-l}(az')],-k-l)              \\
& = & (xz,h^{-(k+l)}(az'),-(k+l))   \\
& = & F((xz,k,yz)(y_{1}z',l,az')).
\eeqns
The other facts to be checked are left to the reader.

This deals with the case where every vertex $v$ is $r(e)$ for some $e\in E$.
For the general case, one considers a larger graph $\calE'$ obtained by
adding an infinite left tail to each $v\in V$ for which
$r^{-1}(\{v\})=\emptyset$, applies the preceding result to $\calE'$ and
performs a reduction to obtain that $G$ is amenable.
(Note that $G^{0}=G(\calE)^{0}$ is a closed subset of $G(\calE')^{0}$.)
We therefore have:

\begin{theorem}                 \label{th:amen}
For any directed graph $\calE$, the graph groupoid $G$ is amenable.
\end{theorem}

\section{The simplicity of $C^{*}(\calE)$}

In this section we will use the groupoid $G$ to obtain conditions sufficient
for $C^{*}(\calE)$ to be simple. As discussed in the Introduction, the result
obtained (Theorem~\ref{th:simple}) is effectively due to W. Szymanski
who also proved the converse.  (It seems likely that the
approach to the converse for the row finite case in \cite{BPRS} can be adapted
to apply within the groupoid context of the present paper.)

For $v\in V$, a {\em loop based at $v$} is a finite path $\al=e_{1}\ldots
e_{n}$ such that $r(\al)=s(\al)=v$ and $r(e_{i})\neq v$ for $1\leq i<l(\al)$.
The subsets $V_{0},V_{1}, V_{2}$ of $V$ are the sets of vertices $v$ such that
there are respectively no loops based at $v$, exactly one loop based at $v$ and
at least two distinct loops based at $v$. The graph $\calE$ is said to
satisfy {\em condition (K)} (\cite[\S 6]{KPRR}) if $V_{1}=\emptyset$.

Recall that a locally compact groupoid $H$ is called {\em essentially
principal} (\cite[p.100]{rg}) if whenever $F$ is a closed invariant subset of
$H^{0}$, then the set of $u$'s in $F$ with trivial isotropy group is dense in
$F$. It is shown in \cite[Proposition 6.3]{KPRR} that if $\calE$ is row
finite and satisfies condition (K), then $G$ is essentially principal.  We
now show that for general $\calE$, $G$ is essentially principal if and only
if $\calE$ satisfies condition (K).

\begin{proposition}
The r-discrete groupoid $G$ is essentially principal if and only if $\calE$
satisfies condition (K).
\end{proposition}
\begin{proof}
Suppose that $\calE$ satisfies condition (K). Let $F$ be a closed invariant
subset of $G$. Note that any finite path $y\in Y$ has trivial isotropy. So we
need only consider infinite paths. The proof of \cite[Proposition 6.3]{KPRR}
then goes through {\em verbatim} to give that $G$ is essentially principal.

Conversely, suppose that $G$ is essentially principal and that $V_{1}\neq
\emptyset$.  Let $v\in V_{1}$ and $\al$ be the loop in $\calE$ based at $v$.
Let
\[       C=\{z\in Z: z\geq v\}.                   \]
Note that $\ga=\al\al\ldots$ belongs to $C$. Let $F=\ov{C}$, the closure of $C$
in $X$. We claim that $F$ is a closed invariant subset of $X$.

Trivially, $F$ is closed. To prove invariance, let $\{z^{i}\}$ be a sequence
in $C$ and $z^{i}\to x\in X$. Suppose that $x$ is equivalent to $y\in Y$.
Either $l(x)=\infty$ or $l(x)<\infty$. Suppose first that $l(x)=\infty$. Then
for some $\al'$ and some $z\in Z$, we have $x=\al' z, y=\bt z$. Then
eventually, every $z^{i}=\al' w^{i}$, and so $\bt w^{i}\in C$ and $\bt w^{i}\to
y$. So $y\in F$. If $l(x)<\infty$, then $r(x)=r(y)$, and a similar argument
gives that $y\in F$. So $F$ is invariant. Note that if $x\in F$, then $x\geq v$
for all $v$ (so that $F\cap Z=C$).

We will contradict the assumption $V_{1}\neq \emptyset$ by showing that if
$x\in F$ and $s(x)=v$, then $x=\ga$. (For then any sequence in $F$ converging
to $\ga$ will eventually coincide with $\ga$ and so eventually,
the terms of the sequence will not have trivial isotropy.)
We can suppose that $x\in C$ since each $x'\in F$ of finite length and with
$s(x')=v$ is the limit of a sequence of such infinite paths $x$. For each $n$,
there is a path $y$ starting and finishing at $v$ obtained by following $x$
until we reach $r(x_{n})$ and then taking a path from $r(x_{n})$ to $v$.
Since $\al$ is the only loop based at $v$, $y$ is of the form $\al\al\ldots
\al$. It follows every initial segment of $x$ is an initial segment of $\ga$
and so $x=\ga$.
\end{proof}

The groupoid $G$ is called {\em minimal} (\cite[p.35]{rg}) if the only open
invariant subsets of $G^{0}$ are $G^{0}$ and the empty set. Since $G$ is
amenable, we have $C^{*}(G)=C^{*}_{red}(G)$ (\cite[p.92]{rg}). Assume that
$\calE$ satisfies condition (K). Then by the above, the r-discrete groupoid $G$
is essentially principal. So by \cite[Proposition 2.4.6, p.103]{rg}, we have
that $C^{*}(G)$ is simple if and only if $G$ is minimal. We now give a groupoid
proof of the simplicity result for $C^{*}(\calE)$ (\cite[Theorem 12]{Szy1}).

\begin{theorem}     \label{th:simple}
The graph $\cs$ $C^{*}(\calE)$ is simple if (a), (b) and (c) hold, where:
\bi
\item[(a)] $\calE$ has property (K);
\item[(b)]  $\calE$ is {\em cofinal} (\cite{KPRR}) in the sense that given
$v\in V$ and $z\in Z$, there exists an $n$ such that $v\geq r(z_{n})$;
\item[(c)] if $v\in V$ and $w\in V_{\infty}$, then there is a path from $v$ to
$w$.
\ei
\end{theorem}
\begin{proof}
Suppose that $G$ satisfies (a), (b) and (c). By the preceding comments, we just
need to show that $G$ is minimal.

Let $U\neq \emptyset$ be an open invariant subset of $X$. By
Proposition~\ref{prop:dense}, $U\cap Z\neq\emptyset$. By considering a
neighborhood in $X$ of some $z\in U\cap Z$ and using (ii) of
Proposition~\ref{prop:unit}, we have that for some $\al$, $D_{\al}'\subset U$.
Let $z\in Z$. By (b), there is a path $\bt$ from $r(\al)$ to $r(z_{n})$ for
some $n$. Then $z'=\al\bt z_{n+1}\ldots\in D_{\al}'\subset U$. Since $z'$ is
tail equivalent to $z$ and $U$ is invariant, we have $z\in U$. Next let $y\in
Y_{\infty}$. By (c), there is a path $\bt'$ from $r(\al)$ to $r(y)$. Then
$\al\bt'$ is equivalent to $y$, and $\al\bt'\in U$. So $y\in U$. So $U=Z\cup
Y_{\infty}=X$ and $G$ is minimal.
\end{proof}


\end{document}